\documentclass[a4paper,11pt]{article}

\usepackage{diagrams}
\usepackage{latexsym}
\usepackage{amsfonts}
\usepackage{amssymb}
\usepackage{amsmath}
\usepackage{amscd}
\usepackage{eucal}
\usepackage{amsthm}
\usepackage[all, knot]{xy}
\xyoption{arc}
\xyoption{rotate}
\xyoption{dvips}
\usepackage[dvips]{pict2e}
\usepackage{verbatim} 
\usepackage[dvips]{graphicx}

\newcommand{\op}{\ensuremath{\mbox{\hspace{1pt}{\scriptsize\upshape op}}}}
\newcommand{\set}{\ensuremath{\mbox{\bfseries Set}}}

\newcommand{\cat}[1]{\ensuremath{\mbox{\bfseries {\upshape {#1}}}}}

\newcommand{\bb}[1]{\ensuremath{\mathbb {#1}}}

\newcommand{\lra}{\ensuremath{\longrightarrow}}
\newcommand{\map}[1]{\ensuremath{\stackrel{{#1}}{\lra}}}

\newtheorem{theorem}{Theorem}[section]

\newtheorem{proposition}[theorem]{Proposition}

\newtheorem{lemma}[theorem]{Lemma}
\newtheorem{definition}[theorem]{Definition}

\newsavebox{\quotename}
\newsavebox{\quoteref}

\newcommand{\numroman}{\renewcommand{\labelenumi}{\roman{enumi})}}
\newcommand{\numarabic}{\renewcommand{\labelenumi}{\arabic{enumi})}}

\newcommand{\pica}{\begin{center} \input}
\newcommand{\picz}{\end{center}}
\newcommand{\length}[1]{\setlength{\unitlength}{#1}}

\newlength{\leng}

\newcommand{\sunit}{\setlength{\unitlength}{1mm}}

\newcommand{\gset}{\cat{GSet}}
\newcommand{\coll}{\cat{Coll}}
\newcommand{\contr}{\cat{Contr}}
\newcommand{\opd}{\cat{Opd}}
\newcommand{\owc}{\cat{OWC}}
\newcommand{\cwc}{\contr}
\newcommand{\contrk}[1]{\ensuremath{\contr_{#1}}}
\newcommand{\opdk}[1]{\ensuremath{\opd_{#1}}}
\newcommand{\cwck}[1]{\ensuremath{\cwc_{#1}}}
\newcommand{\owcij}[2]{\ensuremath{\owc_{#1,#2}}}

\newarrow{Dashto}{}{dash}{}{dash}>
\newarrow{Pointto}----{triangle}
\newarrow{Dashpointto}{}{dash}{}{dash}{triangle}
\newarrow{Into}C--->
\newarrow{Monic}>--->
\newarrow{Implies}===={=>}

\newcommand{\diagsmall}{\begin{diagram}[labelstyle=\scriptstyle]}

\newcommand{\myvec}[3]{\ensuremath{\left( \begin{diagram}[labelstyle=\scriptstyle,height=16pt] {#1} \\ \dTo>{#3} \\ {#2} \end{diagram} \right)}}

\newcommand{\pb}[8]{\ensuremath{
\begin{diagram}[labelstyle=\scriptstyle]
{#1}\SEpbk & \rTo^{#2} & {#3} \\
\dTo<{#8} & & \dTo>{#4} \\
{#7} & \rTo_{#6} & {#5} \\
\end{diagram}}}

\newcommand{\adjn}[4]{\ensuremath{
\diagsmall {#1} & \pile{\rTo^{#2} \\ {\scriptstyle \top} \\ \lTo_{#4}} & {#3}}\end{diagram}}

\newcommand{\psinv}[4]{\ensuremath{
\diagsmall {#1} & \pile{\rTo^{#2} \\ \lTo_{#4}} & {#3}}\end{diagram}}

\usepackage{amssymb}

\begin{document}

\title{Monad interleaving: a construction of the operad for Leinster's weak $\omega$-categories}
\author{Eugenia Cheng\\ \\Department of Pure Mathematics,
University
of Sheffield\\E-mail: e.cheng@sheffield.ac.uk}
\date{15th July, 2008}
\maketitle

\begin{abstract}
We show how to ``interleave'' the monad for operads and the monad for contractions on the category \coll\ of collections, to construct the monad for the operads-with-contraction of Leinster.  We first decompose the adjunction for operads and the adjunction for contractions into a chain of adjunctions each of which acts on only one dimension of the underlying globular sets at a time.  We then exhibit mutual stability conditions that enable us to alternate the dimension-by-dimension free functors.  Hence we give an explicit construction of a left adjoint for the forgetful functor $\owc \lra \coll$, from the category of operads-with-contraction to the category of collections.  By applying this to the initial (empty) collection, we obtain explicitly an initial operad-with-contraction, whose algebras are by definition the weak $\omega$-categories of Leinster.
\end{abstract}

\setcounter{tocdepth}{1}
\tableofcontents

\section*{Introduction}
\addcontentsline{toc}{section}{Introduction}

The aim of this work is to give a dimension-by-dimension construction of the operad for Leinster's weak $\omega$-categories.  This operad was introduced in \cite{lei1} (see also \cite{lei7,lei8}) as the initial object in a certain category of ``operads with contraction''.  The existence of such an initial object is given by abstract considerations but a construction is desirable for the purposes of calculations in the resulting theory of  $\omega$-categories, as well as for comparisons with other theories, especially the related theories of Batanin \cite{bat1} and Penon \cite{pen1}.

Leinster's theory of $\omega$-categories is a ``globular'' one, that is, the underlying data for an $\omega$-category is taken to be a \emph{globular set}; each cell has precisely one source cell and one target cell, whose source and target must match according to the so-called globularity conditions.  The composition and coherence is controlled by a particular ``globular operad''.  Globular operads were introduced by Batanin in \cite{bat1}; they fit into a wider picture of generalised operads introduced by Burroni in \cite{bur1} and described in detail by Leinster in \cite{lei8}.  Globular operads are a generalisation of classical operads, introduced for the purposes of studying higher-dimensional algebra.  Where classical operads govern the weakly associative composition of loops in a space, globular operads govern the weakly associative composition of cells of all dimension in an $\omega$-category.

Batanin's idea in \cite{bat1} is to introduce a class of operads that should ``detect'' weak $\omega$-categories, that is, such that weak $\omega$-categories are precisely the algebras for any of these operads.  The class in question is the class of ``contractible operads with a system of compositions''.  Leinster streamlined this notion by generalising the notion of contraction, subsuming the notion of system of compositions.  Thus we are interested in the category of \emph{operads-with-contraction}.

An operad-with-contraction is, as the name suggests, an operad equipped with the structure of a ``contraction''.  This notion of contraction is similar to the notion of contraction in topology; we will give more introductory explanation of these notions at the beginning of Section~\ref{prelim}, followed by full definitions of all the structures involved.

The important idea for this work is that the notions of ``operad'' and ``contraction'' exist independently -- both have as their underlying data a \emph{collection} (see Section~\ref{colldef}).  In fact, we will define categories
\begin{itemize}
\item \coll\ of collections,
\item \cat{Opd} of globular operads,
\item \cat{Contr} of contractions, and
\item \cat{OWC} of operads-with contraction
\end{itemize}
which fit into the following (strict) pullback square:
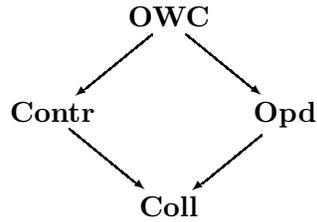
\[\begin{diagram}[labelstyle=\scriptstyle,h=18pt,w=22pt]
&& \cat{OWC} && \\
& \ldTo(2,2) && \rdTo(2,2) & \\
\cat{Contr} &&  && \cat{Opd} \\
& \rdTo(2,2) && \ldTo(2,2) & \\
&& \cat{Coll} &&
\end{diagram}\]
in \cat{CAT} \cite{lei8}.  In particular, an operad-with-contraction is precisely a collection equipped with the structure of both an operad and a contraction, with no further axioms governing the interaction of the two types of structure.

Moreover, it follows from abstract considerations that all the forgetful functors to \coll\ are monadic, including the composite forgetful functor on the ``diagonal''
\[ \cat{OWC} \lra \coll\]
(see \cite[Appendix G]{lei8}).  Our aim is to give convenient constructions of the left adjoints, or equivalently, of the associated monads.  In particular this gives us a construction of the initial object in \cat{OWC}, as we can apply the free functor
\[ \coll \lra \cat{OWC}\]
to the initial object in \coll.  The presence of an initial object in \cat{OWC} means that in practice we do not need to use the whole class of operads to detect $\omega$-categories -- an algebra for any other operad-with-contraction will also be an algebra for the initial one.  Leinster defines $\omega$-categories to be precisely the algebras for this initial operad-with-contraction.

Intuitively, to construct a ``free operad-with-contraction'' monad, we need to start with a collection, and add in both operad and contraction structure freely, using the monad for operads and the monad for contractions.  However, we cannot simply apply one monad and then the other, as we do not have a \emph{distributive law} governing their interaction -- we are simply taking a product of the two monads.  (Note that some authors regard this as a coproduct, e.g. \cite{hm1}, depending on what notion of monad morphism is being used.  We follow Street, as in \cite{str1}.)

The product of these monads may be formed using Kelly's transfinite machinery \cite{kel4}.  However, the aim of this work is to give a more convenient and intuitive construction, hinted at in \cite{lei8} -- we proceed one dimension at a time.  That is, we start by adding in the 0-cells needed for an operad structure, then the 1-cells needed for a contraction structure, then the 1-cells needed for an operad structure, then the 2-cells for a contraction structure, then the 2-cells for an operad structure, and so on.

This raises two questions.

\begin{enumerate}
\item Technical question: why does this dimension-by-dimension construction work?
\item Ideological question: in what way is this construction convenient?
\end{enumerate}
The first question is answered by Propositions~\ref{propone} and \ref{proptwo}.  The idea is that the $k$-cells of a free operad depend only on the lower-dimensional cells, and similarly for a free contraction, so that we can ``decompose'' each individual adjunction into a chain of adjunctions that proceeds one dimension at a time.  We can then alternate or ``interleave'' these two chains of adjunctions, thanks to Lemmas~\ref{lemmaone} and \ref{lemmatwo}, which tell us that the the successive adjunctions of one structure do not interfere with the earlier constructions of the other.

We will now answer the second question somewhat intuitively.  The idea is that the general construction given in \cite{kel4} would require us to repeatedly add in $k$-dimensional structure \emph{even when all the required $k$-dimensional structure is already present}; we would then quotient it back out again.  Our construction makes use of the fact that the $k$-dimensional structure actually remains stable after a certain point in this inductive process.

\subsubsection*{Remarks on Batanin's original definition}

Note that in \cite{bat1} a construction is given of an operad for Batanin's original $\omega$-categories.  Leinster's definition differs from Batanin's in a number of subtle ways, one of which is that in Leinster's variant, contraction cells are \emph{specified}, whereas in Batanin's definition only their \emph{existence} is demanded.  Leinster's definition is concerned with the category of ``operads with specified contraction'' where Batanin's is concerned with the category of ``operads for which a contraction exists but is not actually specified''.  Thus where Leinster seeks an initial object in the category in question, Batanin (necessarily) only asks for a {\em weakly} initial object -- there is a morphism to every other object in the category but it is not unique, so although the operad structure is canonical, the contraction structure is not.

Batanin's construction also differs from ours in that all dimensions are constructed at once.  One consequence of this is that part of the calculation involves constructing coequalisers of operads, a difficult process which can largely be avoided by proceeding dimension-by-dimension.

\subsubsection*{Remarks on generalisation}

\begin{enumerate}

\item Abstractly, we are interested in the product of the two monads in question, in the category of monads on \coll.  This is related to the universal fibre construction of Steenrod \cite{ste1}, which may also be given abstractly as the product of two monads \cite{ber1,bh1}, although in this case the monads are better behaved and hence the construction is simpler.  Steenrod's construction can be generalised to include less well-behaved monads, namely those that do not preserve pushouts; the resulting construction resembles that of Batanin \cite{bat1} rather than the dimension-by-dimension construction given here.

\item Although we give only one specific example of monad interleaving in this work, it is clear how to apply the method to Penon's definition of $\omega$-category \cite{pen1}, interleaving the monads for ``magmas'' and contractions.  This should facilitate a comparison between the two approaches, which we hope to pursue in the future.  Other categories which might be candidates for such constructions include any categories with sets of $k$-cells, such as simplicial sets, opetopic sets, computads and any Reedy categories.

\item An abstract approach to this construction has been given by Hofstra and De Marchi \cite{hm1}, using the framework of indexed monoidal categories and stacks.

\end{enumerate}

\bigskip

The work is structured as follows.  In Section~\ref{prelim} we give the basic definitions of all the structures involved, including some further informal introduction to Leinster's definition of $\omega$-categories.  In Section~\ref{int} we describe the monad interleaving construction in several stages, beginning with an overview of the method used.

Finally we note that this work was first presented at the PSSL79 meeting in Utrecht, 2003, and was posted on the electronic archive in a preliminary form later that year.

\subsubsection*{Acknowledgements}

This work benefited greatly from conversations with Martin Hyland and Tom Leinster; I would also like to thank Michael Batanin for pointing out an interesting error in my original manuscript.

\section{Preliminary definitions} \label{prelim}

In this section we give the basic definitions leading up to and including the operad for Leinster's $\omega$-categories.  We need the notion of ``globular operad'' introduced by Batanin \cite{bat1}, which is a generalisation of the classical notion of operad \cite{bv1,may2}.  The idea is that where a classical operad has operations of arity $k$ for each non-negative integer $k$, a globular operad will have operations of arity $\alpha$ for every \emph{globular pasting diagram} $\alpha$.  A globular pasting diagram is a formal composite of globular cells, which may be depicted by diagrams such as

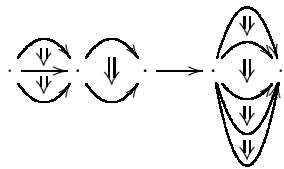
\[\def\objectstyle{\scriptstyle}
       \xy
   (18,0)*+{\cdot}="1";
   (27,0)*+{\cdot}="2";
       {\ar@/^.9pc/ "1";"2"};
       {\ar@/_.9pc/ "1";"2"};
       {\ar@/^2pc/ "1";"2"};
       {\ar@/_2pc/ "1";"2"};
       {\ar@/_3pc/ "1";"2"};
       {\ar@{=>} (22.5,1.5)*{};(22.5,-1.5)*{}} ;
       {\ar@{=>} (22.5,7.25)*{};(22.5,4.75)*{}} ;
       {\ar@{=>} (22.5,-4.75)*{};(22.5,-7.25)*{}} ;
       {\ar@{=>} (22.5,-9.25)*{};(22.5,-11.75)*{}} ;
       (-9,0)*+{\cdot}="1";
       (0,0)*+{\cdot}="2";
           {\ar "1";"2"};
           {\ar@/^1pc/ "1";"2"};
           {\ar@/_1pc/ "1";"2"};
           {\ar@{=>} (-4.5,3)*{};(-4.5,.75)*{}} ;
           {\ar@{=>} (-4.5,-.75)*{};(-4.5,-3)*{}} ;
   (0,0)*+{\cdot}="1";
   (9,0)*+{\cdot}="2";
       {\ar@/^1pc/ "1";"2"};
       {\ar@/_1pc/ "1";"2"};
       {\ar@{=>} (4.5,1.75)*{};(4.5,-1.75)*{}} ;
   (9,0)*+{\cdot}="1";
   (18,0)*+{\cdot}="2";
    {\ar "1";"2"};
\endxy\]
The totality of globular pasting diagrams can be constructed by applying the free strict $\omega$-category monad to the terminal globular set, which has precisely one cell of each dimension.  This is the higher-dimensional version of the fact that the free category monad constructs formal composites of arrows, such as
\[\xy
(0,0)*+{a_0}="1";
(15,0)*+{a_1}="2";
(28,0)*+{}="3";
(35,0)*+{\ldots}="9";
(42,0)*+{}="33";
(57,0)*+{a_{k-1}}="4";
(73,0)*+{a_{k}}="5";
{\ar^{f_1} "1";"2"};
{\ar^{f_2} "2";"3"};
{\ar^{} "33";"4"};
{\ar^{f_k} "4";"5"};
\endxy.\]

The reason for using the generalised notion of operad is that it gives us a convenient way of keeping track of weakly associative composition.  In a strict $\omega$-category, each globular pasting diagram of cells determines precisely one composite, thanks to the strict associativity, unit and interchange axioms.  However in a weak $\omega$-category, we may have many different composites of any given pasting diagram of cells, depending on the order in which the composition is performed.  The different ways of composing a given pasting diagram $\alpha$ can be thought of as ``operations'' of arity $\alpha$, and they form a globular operad; the composition of the operad corresponds to the substitution of one composition scheme into another.

The final subtlety is that the different composites of a given pasting diagram should be related by some notion of equivalence; this is the issue of \emph{coherence} for weak $\omega$-categories.  So we need the notion of ``operad with contraction''; contractibility here is similar to the notion of contractibility of a topological space.  The contraction ensures that there is a coherence cell from every composite of a given diagram to every other; higher-dimensional contractions ensure that every coherence cell is a weak equivalence, as well as ensuring that all coherence cells interact well with each other.

The operad for Leinster's $\omega$-categories is then the initial operad-with-contraction.

\subsection{Globular sets}

The underlying data for an $\omega$-category in this theory is a \emph{globular set}.  This means that each $k$-cell has precisely one $(k-1)$-cell as its source and one as its target, and these in turn have to have source and target matching up according to the \emph{globularity} conditions.  Other definitions of $\omega$-category use more complicated shapes of cells as their underlying data, for example simplicial cells \cite{tam1,str2}, opetopic cells \cite{bd1,hmp3,hmp4,che7} or cubical cells \cite{abs1,bh1,pao1}.

Like simplicial sets, opetopic sets and cubical sets, globular sets are given as presheaves on a category of ``shapes''; in particular this means that the category of globular sets is well-behaved -- it is locally finitely presentable \cite[Example 5.2.2(b)]{bor2}, and so are all the categories derived from it that we use in this work.

\begin{definition}

Let \bb{G} be the category whose objects are the natural numbers $0,1, \ldots$ and whose arrows are generated by
	\[s_k, t_k : k-1 \lra k\]
for each $k \geq 1$, subject to the \emph{globularity equations}
	\[s_k s_{k-1} = t_k s_{k-1} , \ \ s_k t_{k-1}  = t_k t_{k-1} \]
$(k \geq 2)$.  A {\em globular set} is a functor $A: \bb{G}^{\op} \lra \set$, and we write \cat{GSet} for the category of globular sets $[\bb{G}^{\op}, \set]$.
\end{definition}

Explicitly a globular set $A$ consists of a set $A_k$ for each $k \geq 0$, together with morphisms
\begin{diagram}[labelstyle=\scriptstyle,h=28pt,w=28pt]
\cdots  & \pile{\rTo^s \\ \rTo_t} & A_{k+1} & \pile{\rTo^s \\ \rTo_t} & A_k & \pile{\rTo^s \\ \rTo_t} & A_{k-1} & \pile{\rTo^s \\ \rTo_t} & \cdots & \pile{\rTo^s \\ \rTo_t} & A_1 & \pile{\rTo^s \\ \rTo_t} & A_0 \\
\end{diagram}
satisfying $ss=st$ and $ts=tt$.  We refer to the elements of $A_k$ as the \emph{$k$-cells} of $A$. A morphism of globular sets is then a diagram
\begin{diagram}[labelstyle=\scriptstyle,h=28pt,w=28pt]
\cdots  & \pile{\rTo^s \\ \rTo_t} & A_{k+1} & \pile{\rTo^s \\ \rTo_t} & A_k & \pile{\rTo^s \\ \rTo_t} & A_{k-1} & \pile{\rTo^s \\ \rTo_t} & \cdots & \pile{\rTo^s \\ \rTo_t} & A_1 & \pile{\rTo^s \\ \rTo_t} & A_0 \\
\cdots  &  & \dTo>{f_{k+1}} & & \dTo>{f_k} &  & \dTo>{f_{k-1}} & & \cdots & & \dTo>{f_1} &  & \dTo>{f_0} \\
\cdots  & \pile{\rTo^s \\ \rTo_t} & B_{k+1} & \pile{\rTo^s \\ \rTo_t} & B_k & \pile{\rTo^s \\ \rTo_t} & B_{k-1} & \pile{\rTo^s \\ \rTo_t} & \cdots & \pile{\rTo^s \\ \rTo_t} & B_1 & \pile{\rTo^s \\ \rTo_t} & B_0 \\
\end{diagram}
serially commuting.


\subsection{Collections} \label{colldef}

The underlying data for a globular operad is a ``collection'' \cite{bat1}. The idea is to start with a globular set indexed over globular pasting diagrams, constructed using the free strict $\omega$-category monad $T$ on \cat{GSet} (see for example \cite[Appendix G]{lei8} for a construction of this monad).

\begin{definition}
The category \coll\ of {\em collections} is the slice category $\cat{GSet}/T1$.  Here 1 is the terminal globular set, which has precisely one cell of each dimension, and $T$ is the free strict $\omega$-category monad on \gset.
\end{definition}

So a collection is a globular set $A$ together with a morphism $A \map{d} T1$ of globular sets, and a morphism of collections is thus a commuting triangle
%
%
\[\begin{diagram}[labelstyle=\scriptstyle,h=18pt,w=28pt]
A & \rTo^f & B \\
& \rdTo(1,2) \ldTo(1,2) & \\
& T1 & .\\
\end{diagram}\]

\noindent Note that the commutativity is determined dimension-wise.  The $k$-cells of $T1$ are the $k$-dimensional globular pasting diagrams.

The category \coll\ can be given the structure of a monoidal category as follows.  The tensor product of collections $A \map{d} T1$ and $A' \map{d'} T1$ is the composite along the top row of
\begin{diagram}[labelstyle=\scriptstyle,h=28pt,w=28pt]
A \otimes A' \SEpbk & \rTo & TA' & \rTo^{Td'} & T^21 & \rTo^{\mu_1} & T1 \\
\dTo && \dTo>{T!} &&&&\\
A & \rTo^d & T1 &&&&\\
\end{diagram} and the unit for the tensor is $1 \map{\eta_1} T1$ \cite{bat1}.

Note that, as the slice of a presheaf category, \cat{Coll} is itself a presheaf category, and hence also locally finitely presentable \cite{bor2,ar1,lei8}.

\subsection{Globular operads}

Globular operads fit into a wider picture of generalised operads described by Leinster in \cite{lei8}.  This notion of generalised operad was introduced by Burroni in \cite{bur1}; globular operads in particular were introduced by Batanin in \cite{bat1}.

\begin{definition}
A \emph{globular operad} is a monoid in the monoidal category \coll; a {\em morphism of globular operads} is a map of monoids.  We write \cat{Opd} for the category of globular operads and their morphisms, and omit the word ``globular'' for the rest of this work.
\end{definition}

So an operad is given by an underlying collection $A \map{d} T1$ together with a unit and multiplication.  The fibre over a given element $\alpha \in T1$ may be thought of as the set of operations of arity $\alpha$.  For further intuitive explanation of what this data gives us, see \cite{lei8} or \cite{cl1}.  Note that, using the notation of \cite{lei8}, a globular operad is a $(\cat{GSet}, T)$-operad.

\subsection{Contractions}

A contraction is a piece of structure we can look for on any collection $A~\map{d}~T1$.  It is essentially the following lifting property: given a $k$-disc in $T1$, every lift of its boundary to $A$ must give rise to a (specified) lift of the disc itself.  We must now to interpret this in the globular context; this definition first appeared in \cite{lei1}.

\begin{definition}

For any globular set $A$, we call a pair of $k$-cells $a,b \in A_k$ {\em parallel} if $k \geq 1$ and $sa=sb$ and $ta=tb$; all 0-cells are parallel.


Let $A \map{d} T1$ be a collection. Then a {\em contraction} $\gamma$ on it is given by the following data: given

\numroman
\begin{enumerate}
\item a pair of parallel $k$-cells $a,b \in A_k$, and
\item a $(k+1)$-cell $\theta: da \lra db \in T1_{k+1}$
\end{enumerate}
we have a $(k+1)$-cell
	\[\gamma_\theta (a,b) : a \lra b \in A_{k+1}\]
such that $d(\gamma_\theta(a,b)) = \theta$

A {\em collection-with-contraction} is a collection equipped with a specified contraction; for brevity we also refer to a collection-with-contraction simply as a {\em contraction}.  We refer to all the cells $\gamma_\theta(a,b)$ as {\em contraction cells}.

A \emph{morphism of contractions} is a morphism of underlying collections such that the contraction structure is preserved.  Explicitly, given contractions $\gamma$ on $A \map{d} T1$ and $\gamma'$ on $A' \map{d'} T1$, a morphism $\gamma \lra \gamma'$ is a morphism $f: A \lra A'$ of globular sets such that
\begin{enumerate}
\item $d'f=d$, and
\item for all $a,b,\theta$ as above, $f(\gamma_\theta (a,b)) = \gamma'_{\theta}(fa,fb)$.
\end{enumerate}

\noindent We write \cat{Contr} for the category of contractions and their morphisms.
\end{definition}

 Note that Leinster's contractions are more general than Batanin's \cite{bat1} and Penon's \cite{pen1}, which are only required to lift identity cells, rather than all cells $\theta$ as in the above definition.  This is the generalisation that enables Leinster to describe composition and coherence simultaneously using contractions; Batanin and Penon treat composition separately.

\subsection{Operads-with-contraction}

\begin{definition}
An {\em operad-with-contraction} is a collection equipped with both the structure of an operad and the structure of a contraction. A \emph{morphism of operads-with-contraction} is a morphism of underlying collections that is both a morphism of operads and a morphism of contractions.  We write \cat{OWC} for the category of operads-with-contraction and their morphisms.
\end{definition}

Note that we have forgetful functors
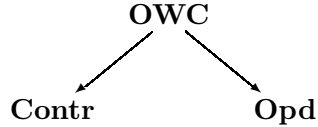
\[\begin{diagram}[labelstyle=\scriptstyle,h=18pt,w=22pt]
&& \cat{OWC} && \\
& \ldTo(2,2) && \rdTo(2,2) & \\
\cat{Contr} &&  && \cat{Opd} \\
\end{diagram}\]

\noindent forgetting just one of the structures at a time.  Forgetting both structures, we get a forgetful functor
	\[\owc \map{G} \coll.\]
Note that \coll\ is locally finitely presentable so it follows from \cite[27.1]{kel4} that $G$ has a left adjoint (see \cite[Appendix G]{lei8}).  We can apply the left adjoint to the initial collection $\emptyset \lra T1$ to obtain an initial object in \owc .

\begin{definition}

The operad for Leinster's $\omega$-categories is the initial operad-with-contraction.

\end{definition}

In this work we will be concerned with giving a convenient construction of the left adjoint to $G$, building it up from the free operad and free contraction functors.  Note that in the definition of operad-with-contraction there are no axioms governing the interaction of these two types of structure; this will later ensure that when we combine the monad for operads and the monad for contractions, we are simply taking a product and not using a distributive law \cite{bec1}.   Put another way, \owc\ is the strict pullback
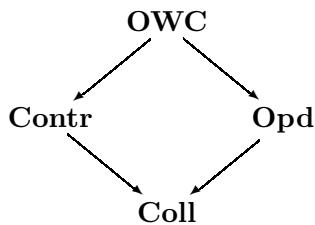
\[\begin{diagram}[labelstyle=\scriptstyle,h=18pt,w=22pt]
&& \cat{OWC} && \\
& \ldTo(2,2) && \rdTo(2,2) & \\
\cat{Contr} &&  && \cat{Opd} \\
& \rdTo(2,2) && \ldTo(2,2) & \\
&& \cat{Coll} &&
\end{diagram}\]
in \cat{CAT} \cite{lei8}, as mentioned in the Introduction.

\subsection{Truncation}

Finally we fix some terminology and notation for the $k$-dimensional versions we will use to build up the $\omega$-dimensional version dimension by dimension.  Note that these are not all the same as the finite-dimensional versions used for defining $n$-categories for finite $n$; the difference is with the contractions, as we will explain below.

\begin{definition}

\ \

\begin{itemize}

\item A {\em $k$-globular set} is a globular set $A$ such that $A_n$ is empty for all $n > k$.

\item A {\em $k$-collection} is a collection whose underlying globular set is $k$-dimensional.

\item A {\em $k$-operad} is an operad whose underlying collection is $k$-dimensional.

\item A {\em $k$-contraction} is a contraction whose underlying collection is $k$-dimensional, but where $k$-cells are not required to have contraction cells between them.


\end{itemize}
\end{definition}
To define $n$-categories for finite $n$, we would need to modify the notion of contraction so that, for $n$-cells, wherever we previously asked for the existence of contraction cells, we now ask for \emph{equalities} between $n$-cells.  It is straightforward to adapt the constructions in this work to the $n$-dimensional version; we will not go into the details.

Note that, when adding structure dimension by dimension, we will have a full (not truncated) underlying collection throughout; only the structure we are adding will be $k$-dimensional at each intermediate stage.  We will use the following notation:

\begin{itemize}

\item For $k \geq 0$ write \opdk{k} for the category of collections whose underlying $k$-collection is equipped with the structure of a $k$-operad.

\item For $k \geq 0$ write \contrk{k} for the category of collections whose underlying $k$-collection is equipped with the structure of a $k$-contraction.

\item For $i,j \geq 0$ write \owcij{i}{j} for the category of collections whose underlying $i$-collection has the structure of an $i$-contraction, and whose underlying $j$-collection has the structure of a $j$-operad.

\end{itemize}  In each case, the morphisms preserve all the structure present.  Note that while $\opdk{0} \neq \coll$, it makes sense to write $\contrk{0}=\coll$.

Finally note that, as complete subcategories of \coll, all the above categories are locally presentable \cite{ar1}, so abstract considerations give the existence of all the adjunctions we use in this work; however, we are interested in constructing the left adjoints in question.

\section{Interleaving structures}\label{int}

In this section we show how to ``interleave'' the operad and contraction structures dimension by dimension, to give a convenient construction of the free operad-with-contraction functor.  As the constructions are somewhat technical we begin with an outline of the method we will adopt.

\subsection{Outline}

\label{outline}

The aim is to combine the two different structures on \coll:

\numroman

\begin{enumerate}
\item operad structure, and
\item contraction structure.
\end{enumerate}
Note that the interaction of these two structures is
different from interactions described by distributive laws;
here the structures have no axioms governing their interaction. This is unlike the notion of a ring, for example, in which we combine the two structures of a monoid and a group, but subject to the distributive law of multiplication over addition, so that the category of rings is not simply a pullback of the categories of monoids and groups.



We will see that we have do two monadic adjunctions

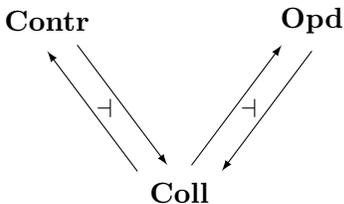
\length{0.4mm}
\begin{center}
\begin{picture}(104,65)(0,-3)

\put(52,1){\makebox(0,0)[c]{$\coll$}}
\put(38,8){\vector(-3,4){30}}
\put(18,50){\vector(3,-4){30}}
\put(27,29){\makebox(0,0)[c]{{\small $\dashv$}}}
\put(8,58){\makebox(0,0)[c]{\contr}}

\put(56,10){\vector(3,4){30}}
\put(96,48){\vector(-3,-4){30}}
\put(75,29){\makebox(0,0)[c]{{\small $\dashv$}}}
\put(96,58){\makebox(0,0)[c]{\opd}}

\end{picture}
\end{center}

\sunit
\noindent and we seek to combine them to get a left adjoint as shown below by the dotted arrow

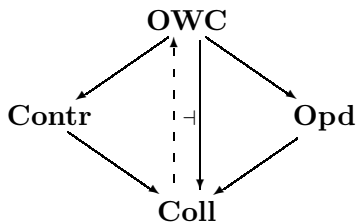
\[\begin{diagram}[labelstyle=\scriptstyle,h=18pt,w=26pt]
&& \cat{OWC} && \\
& \ldTo(2,2) && \rdTo(2,2) & \\
\cat{Contr} && \uDashto>\dashv \dTo && \cat{Opd} \\
& \rdTo(2,2) && \ldTo(2,2) & \\
&& \cat{Coll} && \\
\end{diagram}\]

However, to construct such a left adjoint, we cannot simply proceed `up' one of the sides.  This would amount to adding in the operad structure freely followed by contraction structure freely (or vice versa) and, in effect, the second free structure would destroy the first.  That is, adding in contraction cells \emph{after} the operadic structure would create new operadic composites that would be needed but missing; conversely, adding in the operadic structures after the contraction cells would create new contraction cells that would be needed but missing.

One solution is to alternate the structures transfinitely, but a ``smaller'' solution is available since the new operadic and contraction cells we need to add in are determined only by \emph{lower-dimensional} cells.  This means we can ``interleave'' the structures dimension by dimension, adding in first the free operad structure on 0-cells, then the free contraction structure on 0-cells, then the free operad structure on 1-cells, then the free contraction structure on 1-cells, and so on.

More precisely, we use three key facts:

\numarabic
\begin{enumerate}

\item We have monadic adjunctions
	$\begin{diagram}[labelstyle=\scriptstyle] \cat{Contr} & \pile{\rTo^K \\ {\scriptstyle \top} \\ \lTo_H} &\coll \end{diagram}$
and
	$\begin{diagram}[labelstyle=\scriptstyle] \opd & \pile{\rTo^N \\ {\scriptstyle \top} \\ \lTo_M} &\coll \end{diagram}$
and these restrict to monadic adjunctions on $k$-dimensional truncations.

\item For the free operad on a collection $A$, the new $k$-cells are determined only by the $j$-cells of $A$ for $j \leq k$.

\item For the free contraction on a collection $A$, the new $k$-cells are determined only by the $j$-cells of $A$ for $j \leq k-1$.

\end{enumerate}
The first fact means that each individual construction can proceed dimension by dimension by itself; that is, each adjunction can be decomposed into a chain of adjunctions as follows:

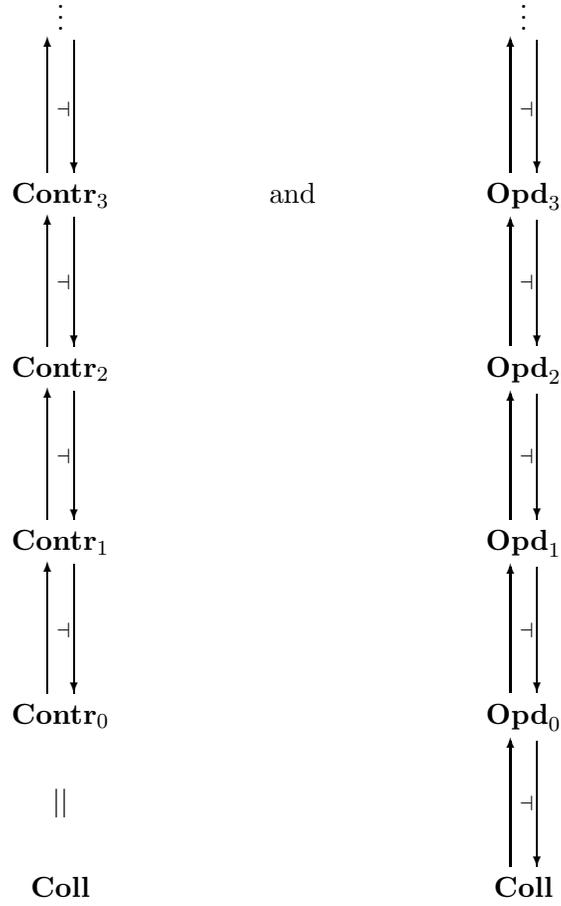
\[\begin{diagram}[labelstyle=\scriptstyle]
\vdots && \vdots \\
\uTo>\dashv \dTo && \uTo>\dashv \dTo \\
\contrk{3} & \mbox{\hspace{2cm} and \hspace{2cm} } & \opdk{3} \\
\uTo>\dashv \dTo && \uTo>\dashv \dTo \\
\contrk{2} && \opdk{2} \\
\uTo>\dashv \dTo && \uTo>\dashv \dTo \\
\contrk{1} && \opdk{1} \\
\uTo>\dashv \dTo && \uTo>\dashv \dTo \\
\contrk{0} && \opdk{0} \\
|| && \uTo>\dashv \dTo \\
\coll && \coll \\
\end{diagram}\]

\bigskip

\noindent The second fact means that a free $k$-operad structure will not be destroyed by adding in new $(k+1)$-cells for a free contraction; the third fact means that a free $k$-contraction structure will not be destroyed by adding in new $k$-cells for a free operad structure.  This makes the interleaving possible, and we can lift the adjunctions as follows:

\[\begin{diagram}[labelstyle=\scriptstyle,h=28pt,w=25pt]
&&&& \vdots &&&& \vdots \\
&&&& \owcij{2}{2} && \rTo && \opdk{2} \\
\vdots &&&& \uTo \dTo<\dashv & {} & \lDashpointto & {} & \uTo \dTo<\dashv \\
\contrk{2} && \lTo && \owcij{2}{1} && \rTo && \opdk{1} \\
\uTo \dTo<\dashv & {} & \rDashpointto & {} & \uTo \dTo<\dashv &&&& ||\\
\contrk{1} && \lTo && \owcij{1}{1} && \rTo && \opdk{1} \\
|| &&&& \uTo \dTo<\dashv & {} & \lDashpointto & {} & \uTo \dTo<\dashv \\
\contrk{1} && \lTo && \owcij{1}{0} && \rTo && \opdk{0} \\
\uTo \dTo<\dashv &{} & \rDashpointto & {} & \uTo \dTo<\dashv &&&& ||\\
\contrk{0} && \lTo && \owcij{0}{0} && \rTo && \opdk{0} \\
&&&& \uTo \dTo<\dashv & {} & \lDashpointto & {} & \uTo \dTo<\dashv \\
&&&& \coll && \rTo^= && \coll \\
&&&& \mbox{Figure 1} &&&& \\
\end{diagram}\]

\bigskip

\noindent Thus the construction proceeds in the following steps.

\numarabic

\begin{enumerate}
\item Dimension-by-dimension decomposition of free operad
construction i.e. the chain of adjunctions on the right hand side above.
\item Dimension-by-dimension decomposition of free contraction
construction i.e. the chain of adjunctions on the left hand side above.
\item Stability of lower-dimensional operad structure under
$k$-dimensional free contraction construction, enabling interleaving of contractions i.e. lifts on left hand side above.
\item Stability of lower-dimensional contraction structure
under $k$-dimensional free operad construction, enabling interleaving of operad structure i.e. lifts on right hand side above.
\item Dimension-by-dimension interleaving of free operad and
free contraction constructions i.e. alternating lifts as in the above diagram.
\end{enumerate}

%

\subsection{Decomposition of the free operad functor}
\label{one}

We know \cite{lei8} that there is a monadic adjunction
	\[\begin{diagram}[labelstyle=\scriptstyle]
	\opd & \pile{\rTo^N \\ {\scriptstyle \top} \\ \lTo_M} & \coll
	\end{diagram}.\]

\begin{proposition}\label{propone}
The above adjunction decomposes as a chain of adjunctions:
\[\begin{diagram}[labelstyle=\scriptstyle]
\cdots & \opdk{k+1} & \pile{\rTo^{N_{k+1}} \\ {\scriptstyle
\top} \\ \lTo_{M_{k+1}}} & \opdk{k} & \pile{\rTo^{N_{k}} \\
{\scriptstyle \top} \\ \lTo_{M_{k}}} & \opdk{k-1} &
\pile{\rTo^{N_{k-1}} \\ {\scriptstyle \top} \\
\lTo_{M_{k-1}}} & \  \cdots \ & \pile{\rTo^{N_{1}} \\ {\scriptstyle
\top} \\ \lTo_{M_{1}}} & \opdk{0} & \pile{\rTo^{N_{0}} \\
{\scriptstyle \top} \\ \lTo_{M_{0}}} & \coll
\end{diagram}\]
with the following properties for each for each $k\geq 0$:

\numroman
\begin{enumerate}
\item $M_k$ leaves all dimensions unchanged except the $k$th dimension, and preserves the underlying $(k-1)$-operad structure, and
\item the underlying $k$-dimensional operad of $MA$ is that of $M_k M_{k-1}\ldots M_0A$.
\end{enumerate}
\end{proposition}

To prove this we copy the construction of the free operad functor
	\[M:\coll \lra \opd\]
\cite{lei8} which builds up formal operadic composites of increasing ``depth'' by induction.  The subtlety here is that at each stage we use the composition of the underlying $k$-operad to keep the underlying $k$-collection stable.  That is, each time we build another depth of formal composite of $(k+1)$-cells, the natural source and target will be a formal composite of $k$-cells, but we can then ``evaluate'' this formal composite, since we already have operadic composition for $k$-cells.

Since an operad is just a monoid in a certain monoidal category, the free operad construction is just a free monoid construction \cite{kel4}; the construction we use here is also just a free monoid construction but in slightly different monoidal category.  Given a $k$-operad $Q$, there is a subcategory of \opdk{k}\ whose objects are those whose underlying $k$-operad is $Q$.  This has a monoidal structure given by using the monoidal structure of \coll\ but quotienting out by the monoid (operad) structure of $Q$ on dimensions $k$ and lower.  We then construct a free monoid in this new monoidal category, which, by construction, leaves the lowest dimensions unchanged as required.

\begin{proof}
First, $M_0(A)$ is given simply by taking the 0-cells of $MA$ and leaving the higher dimensions unchanged; the source and target maps also remain unchanged except at dimension 1, where we must compose with the 0-dimensional component of the unit for the monad $NM$.

Now for each $k \geq 1$ we construct a left adjoint $M_{k+1}: \opdk{k} \lra \opdk{k+1}$, $M_{k+1} \dashv N_{k+1}$.

Let $\myvec{X}{T1}{x} \in \opdk{k}$.  We aim to construct an object
	\[M_{k+1} \myvec{X}{T1}{x} = \myvec{A}{T1}{d} \in \opd_{k+1}\]
where for each $m \neq k+1, A(m) = X(m)$.

We define, for each $n$, a collection $A_n \map{d_n} T1$ of ``composites of depth at most $n$"; we will then construct a chain of inclusions

\[\begin{diagram}[labelstyle=\scriptstyle,h=28pt,w=28pt,heads=littlevee]
A_0 & \rMonic^{i_0} & A_1 & \rMonic^{i_1} & A_2 & \rMonic^{i_2} & \cdots
\end{diagram}\]
and set $A$ to be the colimit of this diagram.

$A_n$ is defined by induction as follows.

\begin{itemize}
\item $A_0 =$

\[\begin{diagram}[labelstyle=\scriptstyle,h=28pt,w=28pt]
\cdots & \pile{\rTo \\ \rTo} & X(k+2) & \pile{\rTo \\ \rTo} & X(k+1) & \rTo^{!} &
1(k+1) & \pile{\rTo \\ \rTo} & 1(k) & \rTo^{u_k} & X(k) & \pile{\rTo \\ \rTo} & \cdots \\
&& \uDashpointto>{} &&&& \uDashpointto>{} &&&& \uDashpointto>{} \\
&& A_0(k+2) &&&& A_0(k+1) &&&& A_0(k) &&
\end{diagram}\]

\noindent Here 1 is the terminal globular set and $u$ is the unit for the underlying $k$-operad of $X$.  We make this into a collection in the obvious way.

\item $A_{n+1}=$

\hspace{-3cm}$\begin{diagram}[labelstyle=\scriptstyle,h=28pt,w=28pt]
\cdots & \pile{\rTo \\ \rTo} & X(k+2) & \pile{\rTo \\ \rTo} & X(k+1) & \rTo^{e_{n+1}} &
(1+X \otimes A_n)(k+1) & \pile{\rTo \\ \rTo} & (1+X \otimes A_n)(k) & \rTo^{\theta_k} & X(k) & \pile{\rTo \\ \rTo} & \cdots \\
&& \uDashpointto>{} &&&& \uDashpointto>{} &&&& \uDashpointto>{} \\
&& A_{n+1}(k+2) &&&& A_{n+1}(k+1) &&&& A_{n+1}(k) &&
\end{diagram}$

\noindent where

\numroman
\begin{enumerate}

\item $e_{n+1}$ is the obvious map given by induction.  Note that eventually this will be used to construct the unit for the adjunction.

\item $\theta_k$ is given by the underlying $k$-operad structure for $X$ since $(1 + X \otimes A_n)(k) = (1 + X \otimes X)(k)$.

\end{enumerate}

\noindent Again, we make this into a collection in the obvious way.

\end{itemize}

\noindent Now define for each $n$ a map $A_n \map{i_n} A_{n+1}$ by the identity at all dimensions not $k+1$, and at the $(k+1)$th dimension:

\begin{itemize}
\item $i_0(k+1): 1(k+1) \lra (1 + X \otimes 1)(k+1)$ is first coprojection
\item $i_{n+1}(k+1): (1+X\otimes A_n)(k+1) \lra (1 + X \otimes A_{n+1})(k+1)$ is $1+X \otimes i_n$.
\end{itemize}
Then the $i_n$'s are monic, and by taking $A$ to be the colimit of
\[\begin{diagram}[labelstyle=\scriptstyle,h=28pt,w=28pt,heads=littlevee]
A_0 & \rMonic^{i_0} & A_1 & \rMonic^{i_1} & A_2 & \rMonic^{i_2} & \cdots
\end{diagram}\]
we obtain a collection \myvec{A}{T1}{}; \coll\ is cocomplete so this limit exists.  Then we can check that this collection naturally has the structure of an operad.

It is then straightforward to construct a unit and counit for the putative adjunction.  The unit comes from the $e_n$ given above.  The counit is also by induction: we seek a morphism
\[\begin{diagram}[labelstyle=\scriptstyle,h=24pt,w=18pt]
A && \rTo^{\epsilon} && X \\
& \rdTo(2,2) && \ldTo(2,2) & \\
&& T1 &&
\end{diagram}\]


\noindent keeping the earlier notation, except that now \myvec{X}{T1}{} already has an underlying  $(k+1)$-operad structure, the $(k+1)$th dimension of which is forgotten for the construction of $A$.  $\epsilon$ will be the indentity at all dimensions except $(k+1)$ where we use the following maps for each $n$:

\begin{itemize}

\item $\epsilon_0:A_0(k+1) = 1(k+1) \map{u} X(k+1)$, the unit for the underlying $(k+1)$-operad of $X$.

\item $\epsilon_{n+1}:\begin{diagram}[labelstyle=\scriptstyle]A_{n+1}(k+1) = (1 + X\otimes A_n)(k+1) & \rTo^{1+X \otimes \epsilon_n} & (1+X \otimes X)(k+1) & \rTo^{\theta_{k+1}} & X(k+1) \end{diagram}$
\end{itemize}

The universal property of the colimit then induces a map $A \map{\epsilon} X$ from the $\epsilon_n$ and we can check that this preserves the underlying $(k+1)$-operad structure.  It is then straightforward to check that the triangle identities are satisfied (pointwise).

 It is immediate from the construction that properties (i) and (ii) hold.
\end{proof}

\subsection{Decomposition of the free contraction functor}
\label{two}

An analogous result holds for contractions.  We first describe the monadic adjunction
	\[\begin{diagram}[labelstyle=\scriptstyle,h=28pt,w=28pt] \cwc & \pile{\rTo^K \\ {\scriptstyle \top} \\ \lTo_H} &
	\coll \end{diagram}\]
\cite{lei8} and then the rest of the construction follows analogously to that of the previous section.

The construction of the adjunction is by induction over dimension.  Given a collection $A \map{f} T1$ we construct a collection-with-contraction $HA \map{\bar{f}} T1$ as follows; the idea is to add in, for each dimension $k$, a set $C_k$ of required contraction cells.

\begin{itemize}

\item $HA_0 = A_0$

\item $HA_1 = A_1 \amalg C_1$ where $C_1$ is given by the pullback
	\[\pb{C_1}{}{HA_0 \times HA_0}{f_0 \times f_0}{T1_0 \times T1_0}{(s,t)}{T1_1}{\bar{f}_1}\]
and we define $(s,t)$ on $C_1$ to be the morphism along the top.  Thus $C_1$ gives all triples $(a,b,\theta)$ requiring a contraction 1-cell $\gamma_\theta(a,b)$.

\item for $k \geq 2$ we have $HA_k=A_k \amalg C_k$ given as follows.  Write $\widetilde{HA}_{k-1}$ for the pullback
	\[\pb{.}{}{HA_{k-1}}{(s,t)}{HA_{k-2} \times HA_{k-2}}{(s,t)}{HA_{k-1}}{}\]
giving all parallel pairs of $(k-1)$-cells.
Then $C_k$ is given by the pullback
	\[\diagsmall
	C_k\SEpbk & \rTo & \widetilde{HA}_{k-1} \\
	&& \dInto \\
	\dTo<{\bar{f}_k} && HA_{k-1} \times HA_{k-1} \\
	&& \dTo>{\bar{f}_{k-1} \times \bar{f}_{k-1}} \\
	T1_k & \rTo_{(s,t)} & T1_{k-1} \times T1_{k-1} \\
\end{diagram}\]
giving all triples $(a,b,\theta)$ requiring a contraction $k$-cell $\gamma_\theta(a,b)$. Then $(s,t)$ on $C_k$ is the composite
	\[C_k \lra \widetilde{HA}_{k-1} \lra HA_{k-1} \times HA_{k-1}.\]
\end{itemize}

Globularity and axioms for a morphism of globular sets follow immediately, and contraction cells are given by the $C_k$.  This clearly gives a monadic adjunction as required, and restricts to a monadic adjunction on $k$-truncations.  Note that each set $C_k$ of contraction $k$-cells is determined only by $j$-cells of $A$ for $j < k$ (together with cells of $T1$).


\begin{proposition}\label{proptwo}
The above adjunction decomposes as a chain of adjunctions
	\[\begin{diagram}[labelstyle=\scriptstyle,h=28pt,w=36pt] \cdots & \cwck{k+1} &
	\pile{\rTo^{K_{k+1}} \\ {\scriptstyle \top} \\ \lTo_{H_{k+1}}} & \cwck{k} &
	\pile{\rTo^{K_{k}} \\ {\scriptstyle \top} \\ \lTo_{H_{k}}} & \cwck{k-1} &
	\pile{\rTo^{K_{k-1}} \\ {\scriptstyle \top} \\ \lTo_{H_{k-1}}} & \mbox{\hspace{2mm}}\cdots \mbox{\hspace{2mm}} &
	\pile{\rTo^{K_{1}} \\ {\scriptstyle \top}\ \ \  \\ \lTo_{H_{1}}} & \cwck{0} = \coll
	\end{diagram}\]	
such that the following two properties hold for each $k\geq 1$:

\numroman
\begin{enumerate}
\item $H_k$ leaves all dimensions unchanged except the $k$th dimension, and preserves the underlying $(k-1)$-contraction structure, and
\item the underlying $k$-dimensional collection-with-contraction of $HA$ is that of $H_k H_{k-1}\ldots H_1A$.
\end{enumerate}
\end{proposition}

\begin{proof}
The construction of each adjunction $H_k \dashv K_k$ proceeds analogously to the constructions given in the previous section.
\end{proof}

\subsection{Interleaving of free contractions}
\label{three}

In this section we give the result about the free contraction functors which makes interleaving possible.  For all $i,j$ we have a forgetful functor
	\[\owcij{i}{j} \lra \cwc_i\]
forgetting the operad structure.  We also have for all $k \geq 0$ a forgetful functor
	\[G_{k+1,k}:\owcij{k+1}{k} \lra \owcij{k}{k}\]
forgetting just the $(k+1)$th-dimension of operad structure.

\begin{lemma}\label{lemmatwo}
For all $k \geq 0$ the adjunction
	\[\begin{diagram}[labelstyle=\scriptstyle] \cwck{k+1} & \pile{\rTo^{K_{k+1}} \\ {\scriptstyle \top} \\
	\lTo_{H_{k+1}}} & \cwck{k} \end{diagram}\]
lifts to an adjunction
	\[\begin{diagram}[labelstyle=\scriptstyle] \owcij{k+1}{k} & \pile{\rTo^{G_{k+1,k}} \\ {\scriptstyle \top} 	
	\\ \lDashto_{F_{k+1,k}}} & \owcij{k}{k} \end{diagram}\]
making the following diagram serially commute:
	\[\begin{diagram}[labelstyle=\scriptstyle]
	\owcij{k+1}{k} & \pile{\rTo^{G_{k+1,k}} \\ {\scriptstyle \top} \\
		\lDashto_{F_{k+1,k}}} & \owcij{k}{k} \\
	\dTo & & \dTo \\
	\cwck{k+1} & \pile{\rTo^{K_{k+1}} \\ {\scriptstyle \top} \\
		\lTo_{H_{k+1}}} & \cwck{k}. \\
	\end{diagram}\]

\end{lemma}

\begin{proof}
We need to show that the $k$-operad structure of any $A \in \cwck{k}$ is ``stable'' under $H_{k+1}$; this is immediate from the fact that, by construction, the underlying $k$-globular set of an $A$ is stable under the action of $H_{k+1}$ as it only adds $(k+1)$-cells.
\end{proof}

%

\subsection{Interleaving of free operad structure}
\label{four}

We now consider the operad structure, and prove the analogous result about the free operad functors which makes interleaving possible.   We have for all $i,j$ a forgetful functor
	\[\owcij{i}{j} \lra \opdk{j}\]
forgetting the contractions, and for all $k\geq 0$ a forgetful functor
	\[G_{k+1,k+1} : \owcij{k+1}{k+1} \lra \owcij{k+1}{k}\]
forgetting just the $(k+1)$the dimension of contraction structure.

\begin{lemma}\label{lemmaone}
For all $k \geq 0$ the adjunction
	\[\begin{diagram}[labelstyle=\scriptstyle] \opdk{k+1} & \pile{\rTo^{N_{k+1}} \\ {\scriptstyle \top} \\
	\lTo_{M_{k+1}}} & \opdk{k} \end{diagram}\]
lifts to an adjunction
	\[\begin{diagram}[labelstyle=\scriptstyle] \owcij{k+1}{k+1} & \pile{\rTo^{G_{k+1,k+1}} \\ {\scriptstyle 	
	\top} \\ \lDashto_{F_{k+1,k+1}}} & \owcij{k+1}{k} \end{diagram}\]
making the following diagram serially commute:
	\[\begin{diagram}[labelstyle=\scriptstyle]
	\owcij{k+1}{k+1} & \pile{\rTo^{G_{k+1,k+1}} \\ {\scriptstyle \top} \\
		\lDashto_{F_{k+1,k+1}}} & \owcij{k+1}{k} \\
	\dTo & & \dTo \\
	\opdk{k+1} & \pile{\rTo^{N_{k+1}} \\ {\scriptstyle \top} \\
		\lTo_{M_{k+1}}} & \opdk{k}. \\
	\end{diagram}\]
\end{lemma}

\begin{proof}
We must show that the $(k+1)$-contraction structure is ``stable'' under $M_{k+1}$.  As above, we know that the underlying $k$-globular set of an object $A \in \opdk{k}$ is stable under the action of $M_{k+1}$; this functor only adds $(k+1)$-cells.  Now suppose that $A$ has the structure of a $(k+1)$-contraction.  Now, since the contraction cells required depend only on the cells of $A$ of dimension $k$ and below (and on $T1$), adding more $(k+1)$-cells does not affect this structure.


So if $A\in \opdk{k}$ has the structure of a $(k+1)$-contraction then so does $M_{k+1}A$, and we have the adjunction as required.
\end{proof}

\subsection{Combining the structures}
\label{five}

We now combine the above results.  We alternate the functors:

\begin{itemize}
\item $F_{k+1,k+1}$ to add operad structure, and
\item $F_{k+1,k}$ to add contraction structure.
\end{itemize}
That is, we have a chain of adjunctions

	\[\begin{diagram}[labelstyle=\scriptstyle]
	\cdots && \owcij{k+1}{k+1} & \pile{\rTo^{G_{k+1,k+1}} \\ {\scriptstyle \top} \\
	\lTo_{F_{k+1,k+1}}} & \owcij{k+1}{k} & \pile{\rTo^{G_{k+1,k}}
	\\ {\scriptstyle \top} 	\\ \lTo_{F_{k+1,k}}} & \owcij{k}{k} & 	
	\pile{\rTo^{G_{k,k}} \\ {\scriptstyle \top} \\ \lTo_{F_{k,k}}} &
	 \owcij{k}{k-1} & \cdots &&&&\\
	&&&&& \cdots & \pile{\rTo^{G_{1,0}} \\ {\scriptstyle \top} \\
	\lTo_{F_{1,0}}} & \owcij{0}{0} & \pile{\rTo^{G_{0,0}} \\ {\scriptstyle \top} 	
	\\ \lTo_{F_{0,0}}} & \coll \end{diagram}\]
given by the central ``spine'' of Figure 1.  So for each $k$ we have a composite adjunction
	\[\begin{diagram}[labelstyle=\scriptstyle]
	\owcij{k}{k} & \pile{\rTo^{\ \ G_k} \\ {\scriptstyle \ \ \top} \\ \lTo_{\ \ F_k}}
	& \coll
	\end{diagram},\]
say. We then define a functor $F:\coll \lra \owc$ as follows.  Let $A$ be a collection.  Put $(FA)_k = (F_kA)_k$ with globular, collection, operad and contraction structures at $k$-dimensions being given by those of $F_kA$.  It then follows that
	\[F \dashv G: \owc \lra \coll\]
as required.

In effect we have contructed a left adjoint for $G$  by taking a limit over the categories \owcij{k}{k}, but the specifics of the construction formalise the intuitive approach previously only hinted at in the literature.

\addcontentsline{toc}{section}{References}


\end{document}